\documentclass[12pt,reqno]{amsart}
\usepackage{enumerate, latexsym, amsmath, amsfonts, amssymb, amsthm, color}
\usepackage{hyperref}
\def\pmod #1{\ ({\rm{mod}}\ #1)}

\def\bg{\bigg}
\def\({\bg(}
\def\){\bg)}

\def\f{\frac}

\def\bi{\binom}

\def\eq{\equiv}

\def\Proof{\noindent{\it Proof}}
\def\Ack{\medskip\noindent {\bf Acknowledgment}}
\theoremstyle{plain}
\newtheorem{theorem}{Theorem}

\newtheorem{lemma}{Lemma}

\newtheorem{conjecture}{Conjecture}
\theoremstyle{definition}

\theoremstyle{remark}
\newtheorem{remark}{Remark}

 \vspace{4mm}

\begin{document}

\hbox{}
\medskip

\title
[{Proof of two conjectures of Z.-W. Sun}]
{Proof of Sun's conjectures on super congruences and the divisibility of certain binomial sums}

\author
[Guo-Shuai Mao and Tao Zhang] {Guo-Shuai Mao and Tao Zhang}

\address {(Guo-Shuai Mao) Department of Mathematics, Nanjing
University, Nanjing 210093, People's Republic of China}
\email{mg1421007@smail.nju.edu.cn}

\address {(Tao Zhang) Department of Mathematics, University of Maryland, College Park,
 MD 20742, USA}
\email{taozhang@math.umd.edu}
\keywords{Central binomial coefficients, congruences, divisibility problem.
\newline \indent This research was supported by the Natural Science Foundation (Grant No. 11571162) of China.}
\begin{abstract}  In this paper, we prove two conjectures of Z.-W. Sun:
$$2n\binom{2n}n\bg|\sum_{k=0}^{n-1}(3k+1)\binom{2k}k^3{16}^{n-1-k}\ \mbox{for \ all}\ n=2,3,\cdots,$$
and
$$\sum_{k=0}^{(p-1)/2}\frac{3k+1}{16^k}\binom{2k}{k}^3\eq p+2\left(\frac{-1}{p}\right)p^3E_{p-3}\pmod{p^4},$$
where $p>3$ is a prime and $E_0,E_1,E_2,\cdots$ are Euler numbers.
\end{abstract}

\maketitle

\section{Introduction}
\setcounter{lemma}{0}
\setcounter{theorem}{0}
\setcounter{corollary}{0}
\setcounter{remark}{0}
\setcounter{equation}{0}
\setcounter{conjecture}{0}
     Let~$p>3$~be a prime. A~$p$-adic congruence is called a super congruence if it happens to hold modulo some higher power of~$p$. Sun \cite{Su1} proved several super congruences involving Euler numbers, such as
\begin{align*}
\sum_{k=0}^{p-1}\frac{\binom{2k}{k}}{2^k}&\eq (-1)^{(p-1)/2}-p^2E_{p-3}\pmod{p^3}.
\end{align*}
Moreover, he proposed many conjectures, such as
\begin{conjecture}\cite[Conjecture 5.1]{Su1} \rm{(i)} For each $n=2,3,\cdots$ we have
 \begin{equation}\label{div1}
   2n\binom{2n}n\bg|\sum_{k=0}^{n-1}(3k+1)\binom{2k}k^3{16}^{n-1-k},
\end{equation}
\begin{equation}\label{div2}
   2n\binom{2n}n\bg|\sum_{k=0}^{n-1}(42k+5)\binom{2k}k^3{4096}^{n-1-k}.
  \end{equation}
 \rm{(ii)} Let $p>3$ be a prime. Then
 \begin{equation}\label{cong1}
   \sum_{k=0}^{(p-1)/2}\frac{3k+1}{16^k}\binom{2k}k^3\eq p+2\left(\frac{-1}{p}\right)p^3E_{p-3}\pmod{p^4},
\end{equation}
\begin{equation}\label{cong2}
   \sum_{k=0}^{p-1}\frac{42k+5}{4096^k}\binom{2k}k^3\eq 5p\left(\frac{-1}{p}\right)-p^3E_{p-3}\pmod{p^4},
 \end{equation}
where~$\left(\frac{\cdot}{p}\right)$~denotes the Legendre symbol.
\end{conjecture}

The congruence conjecture \eqref{cong2} was solved by Hu and the first author \cite{HM}, while the divisibility conjecture \eqref{div2} remains open.

In \cite{Su2}, Z.-W. Sun proved some products and sums divisible by central binomial coefficients, like
$$4(2n+1)\binom{2n}n\bg|\sum_{k=0}^n(4k+1)\binom{2k}k^3(-64)^{n-k}$$
for any positive integer $n$.

Guo also proved some products and sums divisible by central binomial coefficients.
The reader is referred to \cite{G1}.

Motivated by the above work, we obtain the following result.

\begin{theorem}\label{Th1.1} For $n=2,3,4,\ldots$, the assertion \eqref{div1} is true.
 \end{theorem}

Recently, a $q$-analogue of \eqref{div1} has been conjectured by Guo \cite[Conjecture 1.7]{G2}.

Guillera and Zudilin \cite{WZ} proved the weaker version of the congruence conjecture \eqref{cong1}
\[
 \sum_{k=0}^{(p-1)/2}\frac{3k+1}{16^k}\binom{2k}k^3\eq p\pmod{p^3}
\]
using the Wilf-Zeilberger method.

Motivated by their work, we obtain the following result.
\begin{theorem}\label{Th1.2}
 Let~$p>3$~be a prime. Then the congruence \eqref{cong1} is true.
\end{theorem}
We will prove Theorem \ref{Th1.1} and \ref{Th1.2} in Sects. \ref{sect2} and \ref{sect3}, respectively.
\section{Proof of Theorem 1.1}\label{sect2}
 \setcounter{lemma}{0}
\setcounter{theorem}{0}
\setcounter{corollary}{0}
\setcounter{remark}{0}
\setcounter{equation}{0}
\setcounter{conjecture}{0}

Let $p$ be a prime and $n$ a positive integer. Then the $p$-adic evaluation of $n$, denoted by $\operatorname{ord}_p(n)$, is the largest number $s$ such that $p^s|n$.
 \begin{lemma}\label{Lem2.1}
 For any positive integer $n\geq6$ and $n\neq2^m+1$, where $m$ is an integer, we have
 $$n-\operatorname{ord}_2((n-1)!)\geq3,$$
 where $\operatorname{ord}_2((n-1)!)=\sum_{i=1}^{\infty}\lfloor\frac{n-1}{2^i}\rfloor$.
 \end{lemma}
\Proof. Recall the following theorem of Kummer (1852),

\textit{The $p$-adic valuation of the binomial coefficient $\binom{m+n}m$ is equal to the number of 'carry-overs' when performing the addition of $n$ and $m$, written in base $p$.}

Noting that $n\binom{2n}n=\frac{2^n(2n-1)!!}{(n-1)!}$, we get
\[
  \operatorname{ord}_2\left(n\binom{2n}n\right)=n-\operatorname{ord}_2((n-1)!).
\]
Hence, it suffices to show that~$8\mid n\binom{2n}n$.

Write $n=a_1a_2\cdots a_k$~in binary expansion with~$a_1=1$. Since $n\geq6$, we have $k>2$.

\noindent{\it \rm Case 1}. If~$a_k=a_{k-1}=0$, then~$4\mid n$. Since~$2\mid\binom{2n}n$~in any case, we get~$8\mid n\binom{2n}n$.

 \noindent{\it \rm Case 2}. If~$a_k=0, a_{k-1}=1$, then $2\mid n$. Since~$a_1=1$~and~$k>2$, by Kummer's Theorem, we have $4\mid\binom{2n}n$. Therefore,  $8\mid n\binom{2n}n$.\\
 {\it \rm Case 3}. We have $a_k=1$.  Since~$n\neq 2^m+1$ and $k>2$, there must be an integer~$i\in \{2,3,\ldots,k-1\}$ such that $a_i=1$. By Kummer's Theorem, $8\mid n\binom{2n}n$.\qed

\begin{remark}\label{rem}
From Case 3 in the proof, we see that if~$n=2^m+1\geq 6$~for some~$m$, then
\[
n-\operatorname{ord}_2((n-1)!)\geq2.
\]
\end{remark}

\begin{lemma}\label{Lem2.2} Let~$n=2^m+1\geq6$~be an integer, where $m$ is an integer. Then we have
$$8~~\bg|\binom{4n-2}{2n-1}\pm 2\binom{2n-2}{n-1}.$$
\end{lemma}
\Proof.  First, we know
\begin{align*}
\binom{4n-2}{2n-1}=&\frac{(4n-2)(4n-3)\cdots(2n+1)(2n)!}{(2n-1)^2(2n-2)^2\cdots(n+1)^2(n!)^2}=\frac{2^n(4n-3)!!}{(2n-1)!!(n-1)!}
\\=&\frac{2^n(4n-3)(4n-5)\cdots(2n+1)}{(n-1)!},
\end{align*}
hence,
\begin{align*}
&\binom{4n-2}{2n-1}\pm2\binom{2n-2}{n-1}
\\=&\frac{2^n(4n-3)(4n-5)\cdots(2n+1)}{(n-1)!}\pm\frac{2^n(2n-3)!!}{(n-1)!}
\\=&\frac{2^n}{(n-1)!}((4n-3)(4n-5)\cdots(2n+1)\pm(2n-3)!!).
\end{align*}
Noting that $(4n-3)(4n-5)\cdots(2n+1)\pm(2n-3)!!\equiv0\pmod 2$, we can deduce that
$$\operatorname{ord}_2\left(\binom{4n-2}{2n-1}\pm2\binom{2n-2}{n-1}\right)\geq n+1-\operatorname{ord}_2((n-1)!)=n+1-\operatorname{ord}_2(n!)$$ since $n=2^m+1$ is an odd integer.
While $n\geq6$ and $n=2^m+1$ we have
\begin{align*}
n+&1-\operatorname{ord}_2(n!)=2^m+2-\sum_{i=1}^{\infty}\left\lfloor\frac{2^m+1}{2^i}\right\rfloor
\\=&2^m+2-(2^{m-1}+2^{m-2}+\cdots+2+1)=2^m+2-(2^m-1)=3.
\end{align*}
Therefore, $\operatorname{ord}_2\left(\binom{4n-2}{2n-1}\pm2\binom{2n-2}{n-1}\right)\geq3$, i.e. $8\mid\binom{4n-2}{2n-1}\pm2\binom{2n-2}{n-1}$. We finish the proof of Lemma \ref{Lem2.2}. \qed


 \noindent{\it \textbf{Proof of Theorem \ref{Th1.1}}}. By \cite[Lemma 3.2]{MS}, for any integer~$n\geq 2$,
\begin{equation*}
 \sum_{k=0}^{n-1}\binom{n-1}{k}^2\binom{x+k}{2n-1}=\f1{(4n-2)\bi{2n-2}{n-1}}\sum_{k=0}^{n-1}(2x-3k)\bi xk^2\bi{2k}k,
 \end{equation*}
set~$x=-\frac{1}{2}$~in this combinatorial identity, we have
\begin{align*}
\sum_{k=0}^{n-1}\binom{n-1}{k}^2\binom{-\frac12+k}{2n-1}&=\f1{(4n-2)\bi{2n-2}{n-1}}\sum_{k=0}^{n-1}(-1-3k)\bi {-\frac12}k^2\bi{2k}k\\
&=\f{-1}{n\bi{2n}{n}}\sum_{k=0}^{n-1}(3k+1)\frac{\bi{2k}k^3}{16^k}.
\end{align*}
It follows that
\begin{align*}
\frac{\sum_{k=0}^{n-1}(3k+1)\binom{2k}k^3{16}^{n-1-k}}{2n\binom{2n}n}&=-\frac{16^{n-1}}2\sum_{k=0}^{n-1}\binom{n-1}{k}^2\binom{-\frac12+k}{2n-1}\\
&=\frac18\sum_{k=0}^{n-1}\binom{n-1}k^2\frac{(-1)^{k}(2k)!(4n-2k-2)!}{k!(2n-1)!(2n-k-1)!}.
\end{align*}
Therefore, to prove Theorem \ref{Th1.1}, we just need to show that
 \begin{equation}\label{eq2.1}
\frac18\sum_{k=0}^{n-1}\binom{n-1}k^2\frac{(-1)^{k}(2k)!(4n-2k-2)!}{k!(2n-1)!(2n-k-1)!}\in\mathbb{Z}.
\end{equation}

When~$n=2,3,4,5$, it is easy to check that \eqref{eq2.1} holds. From now on, we can assume~$n\geq 6$.
For convenience, let
\[
a(n,k)=\frac{(-1)^{k}(2k)!(4n-2k-2)!}{k!(2n-1)!(2n-k-1)!}.
\]
 For any real numbers $x$ and $y$, we have
 $$\lfloor2x\rfloor+\lfloor2y\rfloor\geq \lfloor x\rfloor+\lfloor y\rfloor+\lfloor x+y\rfloor.$$
 It follows that, for any prime~$p$, we have
\begin{align*}
 \operatorname{ord}_{p}\left(a(n,k)\right)=&\sum_{i=1}^{\infty}\left(\left\lfloor\frac{2k}{p^i}\right\rfloor+\left\lfloor\frac{4n-2k-2}{p^i}\right\rfloor-\left\lfloor\frac{k}{p^i}\right\rfloor\right.\\&\left.-\left\lfloor\frac{2n-1}{p^i}\right\rfloor-\left\lfloor\frac{2n-k-1}{p^i}\right\rfloor\right)\\
\geq&0,
\end{align*}
i.e.~$a(n,k)\in\mathbb{Z}$.

Noting that~$$a(n,k)=\frac{(-1)^k2^n(2k-1)!!(4n-2k-3)!!}{(2n-1)!!(n-1)!}.$$
Hence, $$\operatorname{ord}_2\left(a(n,k)\right)=n-\operatorname{ord}_2((n-1)!).$$

If~$n\neq 2^m+1$, by Lemma \ref{Lem2.1} we have $8\mid a(n,k)$.

If~$n=2^m+1$, then for~$1\leq k\leq n-2$, we have $2\mid\binom{n-1}{k}$ and, by Remark \ref{rem}, $4\mid a(n,k)$.

Hence,~$8\mid\binom{n-1}k^2a(n,k)$. For~$k=0$~and~$k=n-1$, note that
\[
\binom{n-1}0^2a(n,0)+\binom{n-1}{n-1}^2a(n,n-1)=\binom{4n-2}{2n-1}\pm 2\binom{2n-2}{n-1},
\]
which is divisible by~$8$~according to Lemma \ref{Lem2.2}.

Therefore,  for any integer $n\geq6$, we have
$$\frac18\sum_{k=0}^{n-1}\binom{n-1}k^2\frac{(-1)^{k+1}(2k)!(4n-2k-2)!}{k!(2n-1)!(2n-k-1)!}\in\mathbb{Z},$$
which completes the proof of Theorem \ref{Th1.1}.\qed


\section{Proof of Theorem 1.2}\label{sect3}
\setcounter{lemma}{0}
\setcounter{theorem}{0}
\setcounter{corollary}{0}
\setcounter{remark}{0}
\setcounter{equation}{0}
\setcounter{conjecture}{0}

\begin{lemma}\label{lem3.1}
Let~$p>3$~be a prime. For any~$0<k\leq(p-1)/2$, we have
\begin{equation}\label{eq3.1}
 \frac{1}{p}\binom{p-1+2k}{(p-1)/2+k}\eq\left(\frac{-1}{p}\right)4^{p-1}
\frac{4^{2k}}{2k\binom{2k}{k}}(1-p(H_{2k-1}-H_{k-1}))\pmod{p^2}.
\end{equation}
In particular,
\begin{equation}\label{eq3.2}
 \frac{1}{p}\binom{p-1+2k}{(p-1)/2+k}\eq\left(\frac{-1}{p}\right)\frac{4^{2k}}{2k\binom{2k}{k}}\pmod{p}.
\end{equation}

\end{lemma}

\begin{proof} Recall that Morley \cite{M} proved that
\begin{equation}\label{Morley}
 \binom{p-1}{(p-1)/2}\eq\left(\frac{-1}{p}\right)4^{p-1}\pmod{p^3}.
\end{equation}
for any prime~$p>3$.
 Hence,
\begin{align*}
 \frac{1}{p}\binom{p-1+2k}{(p-1)/2+k}&=\binom{p-1}{(p-1)/2}\frac{(p+1)\cdots(p+2k-1)}{(\prod_{j=1}^k(p+2j-1)/2)^2}\\
&\eq\left(\frac{-1}{p}\right)4^{p-1}\frac{(2k-1)!(1+pH_{2k-1})2^{2k}}{\prod_{j=1}^{k}(p+2j-1)^2}\\
&\eq\left(\frac{-1}{p}\right)4^{p-1}\frac{(2k-1)!(1+pH_{2k-1})2^{2k}}{\prod_{j=1}^{k}((2j-1)^2+2p(2j-1))}\\
&\eq\left(\frac{-1}{p}\right)4^{p-1}\frac{(2k-1)!(1+pH_{2k-1})2^{2k}}{((2k-1)!!)^2(1+2p\sum_{j=1}^{k}1/(2j-1))}\\
&=\left(\frac{-1}{p}\right)\frac{4^{p-1+2k}(1+pH_{2k-1})}{2k\binom{2k}{k}(1+2p(H_{2k-1}-H_{k-1}/2))}\pmod{p^2}.
\end{align*}
Noting that
\begin{align*}
 \frac{1}{1+2p(H_{2k-1}-H_{k-1}/2)}&=\frac{1-2p(H_{2k-1}-H_{k-1}/2)}{1-4p^2(H_{2k-1}-H_{k-1}/2)^2}\\
&\eq1-2p(H_{2k-1}-H_{k-1}/2)\pmod{p^2}.
\end{align*}
Hence,
\begin{align*}
 \frac{1}{p}\binom{p-1+2k}{(p-1)/2+k}&\eq\left(\frac{-1}{p}\right)\frac{4^{p-1+2k}}{2k\binom{2k}{k}}(1+pH_{2k-1})(1-2p(H_{2k-1}-H_{k-1}/2))\\
&\eq\left(\frac{-1}{p}\right)\frac{4^{p-1+2k}}{2k\binom{2k}{k}}(1-p(H_{2k-1}-H_{k-1}))\pmod{p^2},
\end{align*}
as desired. The last statement follows immediately from Fermat's Little Theorem.
\end{proof}

\begin{lemma}\label{lem3.2}
 Let~$p>3$~be a prime. For any~$0<k\leq(p-1)/2$, we have
\begin{equation}\label{eq3.4}
 \frac{1}{p}\binom{p-1+2k}{2k}\eq \frac{1}{2k}(1+pH_{2k-1})\pmod{p^2}.
\end{equation}
In particular,
\begin{equation}\label{eq3.5}
 \frac{1}{p}\binom{p-1+2k}{2k}\eq \frac{1}{2k}\pmod{p}.
\end{equation}

\end{lemma}
\begin{proof}
 Expanding the LHS, we have
\begin{align*}
 \frac{1}{p}\binom{p-1+2k}{2k}=&\frac{(2k+1)\cdots(p-1)}{1\cdots(p-2k-1)}\cdot\frac{1}{p-2k}\cdot\frac{(p+1)\cdots(p+2k-2)}{(p-2k+1)\cdots(p-1)}\\
&\cdot\frac{1}{2k}\prod_{j=1}^{p-2k-1}(1-p/j)\cdot\frac{1}{1-p/2k}\cdot\prod_{j=1}^{2k-1}\frac{1+p/j}{1-p/j}\\
\eq& \frac{1}{2k}(1-pH_{p-1-2k})(1+p/2k)(1+2pH_{2k-1})\\
\eq& \frac{1}{2k}(1-pH_{p-1-2k}+p/2k+2pH_{2k-1})\pmod{p^2}.
\end{align*}
Recall that Wolstenholem \cite{W} proved that for any prime~$p>3$,
\[
 H_{p-1}\eq0\pmod{p^2}.
\]
If follows that
\begin{align*}
 H_{p-2k-1}&\eq -H_{p-1}+H_{p-2k-1}=-\sum_{j=1}^{2k}1/(p-j)\\
&\eq -\sum_{j=1}^{2k}1/(-j)=H_{2k}\pmod{p}.
\end{align*}
Therefore,
\[
 \frac{1}{p}\binom{p-1+2k}{2k}\eq\frac{1}{2k}(1+pH_{2k-1})\pmod{p^2}.
\]

\end{proof}

\begin{lemma}\label{lem3.3}
 Let~$p>3$~be a prime. Then
\begin{equation}\label{eq3.6}
 \sum_{k=1}^{(p-1)/2}\frac{4^{2k}}{k^2\binom{2k}{k}^2}\eq (-1)^{(p-1)/2}\frac{3}{p}4^{1-p}\sum_{k=1}^{(p-1)/2}\frac{\binom{2k}{k}}{k}\pmod{p}.
\end{equation}
\end{lemma}
\Proof. Noting that
$$\binom{(p-1)/2}{k}\eq\binom{-1/2}{k}=\frac{\binom{2k}{k}}{(-4)^k}\pmod{p},$$
we have
\begin{align*}
 \sum_{k=1}^{(p-1)/2}\frac{4^{2k}}{k^2\binom{2k}{k}^2}&=\sum_{k=1}^{(p-1)/2}\frac{(-4)^{2k}}{k^2\binom{2k}{k}^2}\eq\sum_{k=1}^{(p-1)/2}\frac{1}{k^2\binom{(p-1)/2}{k}^2}\pmod{p}.
\end{align*}
Recall that Staver \cite{S} proved that
\begin{equation*}\label{Staver}
 \sum_{k=1}^{n}\frac{\binom{2k}{k}}{k}=\frac{n+1}{3}\binom{2n+1}{n}\sum_{k=1}^{n}\frac{1}{k^2\binom{n}{k}^2},~\forall n\in\mathbb{Z}^+.
\end{equation*}
Therefore
\begin{align*}
 \sum_{k=1}^{(p-1)/2}\frac{4^{2k}}{k^2\binom{2k}{k}^2}&\eq\frac{3}{\frac{p+1}{2}\binom{p}{(p-1)/2)}}\sum_{k=1}^{(p-1)/2}\frac{\binom{2k}{k}}{k}\\
&=\frac{3}{\frac{p+1}{2}\frac{p}{(p+1)/2}\binom{p-1}{(p-1)/2}}\sum_{k=1}^{(p-1)/2}\frac{\binom{2k}{k}}{k}\\
&\eq(-1)^{(p-1)/2}\frac{3}{p}4^{1-p}\sum_{k=1}^{(p-1)/2}\frac{\binom{2k}{k}}{k}\pmod{p},
\end{align*}
where we use Morley congruence \eqref{Morley} in the last step. \qed

\noindent{\it \textbf{Proof of Theorem \ref{Th1.2}}}.
Take the same WZ pair~$F(k,j)$~and~$G(k,j)$~as in \cite{WZ},
\begin{align*}
 F(k,j)&=\frac{2k+2j+1}{16^k}\binom{2k}{k}^2\frac{\binom{2k+2j}{k+j}\binom{2k+2j}{2j}}{\binom{2j}{j}},\\
 G(k,j)&=-\frac{2(2k-1)}{16^{k-1}}\binom{2k-2}{k-1}^2\frac{\binom{2k+2j-2}{k+j-1}\binom{2k+2j-2}{2j}}{\binom{2j}{j}}.
\end{align*}
We know that~$F(k,j)$~and~$G(k,j)$~have the following relation,
\[
 F(k,j-1)-F(k,j)=G(k+1,j)-G(k,j).
\]
Summing up the above equation for $k$ from~$0$~to~$(p-1)/2$, and then for $j$ from $1$~to~$(p-1)/2$, we get
\[
 \sum_{k=0}^{(p-1)/2}(F(k,0)-F(k,(p-1)/2))=\sum_{j=1}^{(p-1)/2}(G((p+1)/2,j)-G(0,j)).
\]
Noting that~$G(0,j)=0$~and
\[
 \sum_{k=0}^{(p-1)/2}\frac{3k+1}{16^k}\binom{2k}k^3=\sum_{k=0}^{(p-1)/2}F(k,0),
\]
we have
\[
\sum_{k=0}^{(p-1)/2}\frac{3k+1}{16^k}\binom{2k}k^3=\sum_{k=0}^{(p-1)/2}F(k,(p-1)/2)+\sum_{j=1}^{(p-1)/2}G((p+1)/2,j).
\]
Hence, it suffices to determine
\[
 \sum_{k=0}^{(p-1)/2}F(k,(p-1)/2)~\text{and}~\sum_{j=1}^{(p-1)/2}G((p+1)/2,j)\pmod{p^4}.
\]
First, let us consider
\begin{align*}
 G((p+1)/2,j)&=-\frac{2p}{4^{p-1}}\binom{p-1}{(p-1)/2}^2\frac{\binom{p-1+2j}{(p-1)/2+j}\binom{p-1+2j}{2j}}{\binom{2j}{j}}\\
&=-\frac{2p^3}{4^{p-1}}\cdot\frac{\binom{p-1}{(p-1)/2}^2}{\binom{2j}{j}}\cdot\frac{\binom{p-1+2j}{(p-1)/2+j}}{p}\cdot\frac{\binom{p-1+2j}{2j}}{p}.
\end{align*}
By \eqref{eq3.2}, \eqref{Morley}, \eqref{eq3.5} and \eqref{eq3.6}, we get
\begin{align*}
 \sum_{j=1}^{(p-1)/2}G((p+1)/2,j)&\eq-\frac{3}{2}p^2\sum_{j=1}^{(p-1)/2}\frac{\binom{2j}{j}}{j}\pmod{p^4}.
\end{align*}
Now, let us consider
\[
 F(k,(p-1)/2)=\frac{3k+p}{16^k}\binom{2k}{k}^2\frac{\binom{2k+p-1}{k+(p-1)/2}\binom{2k+p-1}{2k}}{\binom{p-1}{(p-1)/2}}.
\]
By \eqref{eq3.1}, \eqref{eq3.2}, \eqref{Morley}, \eqref{eq3.4}, \eqref{eq3.5} and \eqref{eq3.6}, we get
\begin{align*}
&\sum_{k=0}^{(p-1)/2}F\left(k,\frac{p-1}{2}\right)
\\&= p+\sum_{k=1}^{(p-1)/2}\frac{p}{16^k}\binom{2k}{k}^2\frac{\binom{2k+p-1}{k+(p-1)/2}\binom{2k+p-1}{2k}}{\binom{p-1}{(p-1)/2}}
\\
&~~~~+\sum_{k=1}^{(p-1)/2}\frac{3k}{16^k}\binom{2k}{k}^2\frac{\binom{2k+p-1}{k+(p-1)/2}\binom{2k+p-1}{2k}}{\binom{p-1}{(p-1)/2}}\\
&\eq p+\frac{p^3}{4}\sum_{k=1}^{(p-1)/2}\frac{\binom{2k}{k}}{k^2}
+\frac{3p^2}{4}\sum_{k=1}^{(p-1)/2}\frac{\binom{2k}{k}}{k}(1-pH_{2k-1}+pH_{k-1})(1+pH_{2k-1})\\
&\eq p+\frac{p^3}{4}\sum_{k=1}^{(p-1)/2}\frac{\binom{2k}{k}}{k^2}
+\frac{3p^2}{4}\sum_{k=1}^{(p-1)/2}\frac{\binom{2k}{k}}{k}(1+pH_{k-1})\pmod{p^4}.
\end{align*}
Combining them together, we have
\begin{align*}
 \sum_{k=0}^{(p-1)/2}\frac{3k+1}{16^k}\binom{2k}k^3&\eq p-\frac{3}{4}p^2\sum_{k=1}^{(p-1)/2}\frac{\binom{2k}{k}}{k}-\frac{1}{2}p^3\sum_{k=1}^{(p-1)/2}\frac{\binom{2k}{k}}{k^2}\\
&+\frac{3}{4}p^3\sum_{k=1}^{(p-1)/2}\frac{\binom{2k}{k}}{k}H_{k}\pmod{p^4}.
\end{align*}
Finally, applying the congruence (\cite[(1.2)]{Su1})
$$\sum_{k=1}^{(p-1)/2}\frac{\binom{2k}{k}}{k}\eq (-1)^{(p+1)/2}\frac{8}{3}pE_{p-3}\pmod{p^2},$$
and the congruence (\cite[(2.10)]{MS})
$$\sum_{k=1}^{(p-1)/2}\frac{\binom{2k}{k}}{k}H_k\eq \frac{2}{3}\sum_{k=1}^{(p-1)/2}\frac{\binom{2k}{k}}{k^2}\pmod {p}.$$
we get the desired result. \qed

\Ack. The authors would like to thank Professor Z.-W. Sun for helpful comments.

\end{document}